\documentclass[12pt]{article}
\usepackage{hyperref}
\usepackage{amsfonts}
\usepackage{enumerate}
\usepackage{graphicx}
\usepackage{amsmath}
\begin{document}

\begin{center}
{\large\bf The case of a generalised linear discrete time system with infinite many solutions}

\vskip.20in
Nicholas Apostolopoulos$^{1}$, Fernando Ortega$^{2}$ and\ Grigoris Kalogeropoulos$^{3}$\\[2mm]
{\footnotesize
$^{1}$National Technical University of Athens, Greece\\
$^{2}$ Universitat Autonoma de Barcelona, Spain\\
$^{3}$National and Kapodistrian University of Athens, Greece}
\end{center}

{\footnotesize
\noindent
\textbf{Abstract:}
In this article we study a class of generalised linear systems of difference equations with given non-consistent initial conditions and infinite many solutions. We take into consideration the case that the coefficients are square constant matrices with the leading coefficient singular. We provide optimal solutions and numerical examples to justify our theory.
\\\\[3pt]
\\[3pt]

\vskip.2in

\section{Introduction}
Many authors have studied generalised discrete \& continuous time systems, see [1-28], and their applications, see [29-36]. Many of these results have already been extended to systems of differential \& difference equations with fractional operators, see [37-45]. We consider the generalised discrete time system of the form
\begin{equation}\label{eq1}
\begin{array}{cc}
FY_{k+1}=GY_k+V_k, & k= 1, 2,...,
\end{array}
\end{equation}
and the known initial conditions (IC)
\begin{equation}\label{eq2}
Y_{0}.
\end{equation}      
Where  $F, G \in \mathbb{R}^{r \times m}$, $Y_k\in \mathbb{R}^{m}$, and $V_k\in \mathbb{R}^{r}$. The matrices $F$, $G$ can be non-square ($r\neq m$) or square ($r=m$) with $F$ singular (det$F$=0). 

Generalised linear systems of difference equations with given initial conditions don't always guarantee to have unique solution. In the case where there exist solutions and they are infinite, we require optimal solutions for the system. The aim of this paper is to generalise existing results regarding the literature. An explicit and easily testable formula is derived of an optimal solution for the system. 


\section{Preliminaries}
Throughout the paper we will use in several parts matrix pencil theory to establish our results. A matrix pencil is a family of matrices $sF-G$, parametrised by a complex number $s$, see [46-53].
\\\\
\textbf{Definition 2.1.} Given $F,G\in \mathbb{R}^{r \times m}$ and an arbitrary $s\in\mathbb{C}$, the matrix pencil $sF-G$ is called:
\begin{enumerate}
\item Regular when  $r=m$ and  det$(sF-G)\neq 0$;
\item Singular when  $r\neq m$ or  $r=m$ and det$(sF-G)\equiv 0$.
\end{enumerate}
In this article we consider the system \eqref{eq1} with a \textsl{regular pencil}, where the class of $sF-G$ is characterised by a uniquely defined element, known as the Weierstrass canonical form, see [46-53], specified by the complete set of invariants of $sF-G$. This is the set of elementary divisors of type  $(s-a_j)^{p_j}$, called \emph{finite elementary divisors}, where $a_j$ is a finite eigenvalue of algebraic multiplicity $p_j$ ($1\leq j \leq \nu$), and the set of elementary divisors of type $\hat{s}^q=\frac{1}{s^q}$, called \emph{infinite elementary divisors}, where $q$ is the algebraic multiplicity of the infinite eigenvalue. $\sum_{j =1}^\nu p_j  = p$ and $p+q=m$.
\\\\
From the regularity of $sF-G$, there exist non-singular matrices $P$, $Q$ $\in \mathbb{R}^{m \times m}$ such that 
\begin{equation}\label{eq3}
\begin{array}{c}PFQ=\left[\begin{array}{cc} I_p&0_{p,q}\\0_{q,p}&H_q\end{array}\right],
\\\\
PGQ=\left[\begin{array}{cc} J_p&0_{p,q}\\0_{q,p}&I_q\end{array}\right].\end{array}
\end{equation}
$J_p$, $H_q$ are appropriate matrices with $H_q$ a nilpotent matrix with index $q_*$, $J_p$ a Jordan matrix and $p+q=m$. With $0_{q,p}$ we denote the zero matrix of $q\times p$. The matrix $Q$ can be written as
\begin{equation}\label{eq4}
Q=\left[\begin{array}{cc}Q_p & Q_q\end{array}\right].
\end{equation}
$Q_p\in \mathbb{R}^{m \times p}$ and $Q_q\in \mathbb{R}^{m \times q}$. The matrix $P$ can be written as
\begin{equation}\label{eq5}
P=\left[\begin{array}{c}P_1 \\ P_2\end{array}\right].
\end{equation}
$P_1\in \mathbb{R}^{p \times r}$ and $P_2\in \mathbb{R}^{q \times r}$.
The following results have been proved.
\\\\
\textbf{Theorem 2.1.}  (See [1-28]) We consider the systems \eqref{eq1} with a regular pencil. Then, its solution exists and for $k\geq 0$, is given by the formula
\begin{equation}\label{eq6}
    Y_k=Q_pJ_p^kC+QD_k.  
\end{equation}
Where $D_k=\left[
\begin{array}{c} \sum^{k-1}_{i=0}J_p^{k-i-1}P_1V_i\\-\sum^{q_{*}-1}_{i=0}H_q^iP_2V_{k+i}
\end{array}\right]$ and $C\in\mathbb{R}^p$ is a constant vector. The matrices $Q_p$, $Q_q$, $P_1$, $P_2$, $J_p$, $H_q$ are defined by \eqref{eq3}, \eqref{eq4}, \eqref{eq5}. 
\\\\
\textbf{Definition 2.2.} Consider the system \eqref{eq1}  with known IC of type \eqref{eq2}. Then the IC are called consistent if there exists a solution for the system \eqref{eq1} which satisfies the given conditions.
\\\\
\textbf{Proposition 2.1.} The IC \eqref{eq2} of system \eqref{eq1} are consistent if and only if 
\[
Y_0\in colspanQ_p +QD_0.
\]
\textbf{Proposition 2.2.}  Consider the system \eqref{eq1} with given IC. Then the solution for the initial value problem \eqref{eq1}--\eqref{eq2} is unique if and only if the IC are consistent. Then, the unique solution is given by the formula
\[
    Y_k=Q_pJ_p^kZ^p_{0}+QD_k.  
\]
where $D_k=\left[
\begin{array}{c} \sum^{k-1}_{i=0}J_p^{k-i-1}P_1V_i\\-\sum^{q_{*}-1}_{i=0}H_q^iP_2V_{k+i}
\end{array}\right]$ and $Z^p_0$ is the unique solution of the algebraic system $Y_0=Q_pZ^p_0+D_0$.

\section{Non-consistent initial conditions}

We begin this Section with some already existing results.
\\\\
\textbf{Corollary 3.1.} The initial conditions \eqref{eq2} of the system \eqref{eq1} for $V_k=0_{r,1}$ are consistent if and only if 
\[
Y_0\in colspan Q_p.
\]
We can now state the following Theorem. 
\\\\
\textbf{Theorem 3.1.} We consider the system \eqref{eq1} with $V_k=0_{r,1}$ and known non-consistent initial conditions ($Y_0\notin colspanQ_p$) of type \eqref{eq2}. For the case that the pencil $sF-G$ is regular, after a perturbation to the non-consistent initial conditions \eqref{eq2} accordingly
\[
min\left\|Y_0-\hat Y_0\right\|_2,
\]
or, equivalently,
\[
\left\|Y_0-Q_p(Q_p^*Q_p)^{-1}Y_0\right\|_2,
\]
an optimal solution of the initial value problem \eqref{eq1}--\eqref{eq2} is given by
\begin{equation}\label{eq7}
    \hat Y_k=Q_pJ_p^k(Q_p^*Q_p)^{-1}Q_p^*Y_0.
\end{equation}
The matrices $Q_p$, $J_p$ are given by by \eqref{eq3}, \eqref{eq4}, \eqref{eq5}.\\\\
\textbf{Proof.} From Theorem 2.1 the solution of the system \eqref{eq1} exists and is given by \eqref{eq6}. Note that where $D_k=0_{m,1}$, since it is assumed $V_k=0_{r,1}$. For given initial conditions of type \eqref{eq2}, $C$ is the solution of the linear system 
\begin{equation}\label{eq8}
Q_pC=Y_0.
\end{equation}
Since the matrix $Q$ is non-singular, the matrix $Q_p$ has linear independent columns, i.e. $rank Q_p=p$. The system has $m$ linear equations with $p$ unknowns and $m>p$, i.e. the system is overdetermined. To sum up, $Q_p$ is a tall matrix (more rows than columns) with linearly independent columns. Thus, since $Y_0\notin colspanQ_p$ and $rank Q_p=p$, the system \eqref{eq8} has no solutions. Let $\hat Y_0$ be a vector such that $\hat Y_0\in colspanQ_p$ and let $\hat C$ be the unique solution of the system $Q_p\hat C=\hat Y_0$. Hence we want to solve the following optimisation problem
\[
\begin{array}{c}min\left\|Y_0-\hat Y_0\right\|_2^2\\s.t.\quad Q_p\hat C=\hat Y_0,\end{array}
\]
or, equivalently,
\[
min\left\|Y_0-Q_p\hat C\right\|_2^2.
\]
Where $\hat C$ is the optimal solution, in terms of least squares, of the linear system \eqref{eq8}. Thus we seek a solution $\hat C$ by minimising the functional
\[
D_1(\hat C)=\left\|Y_0-Q_p\hat C\right\|_2^2.
\]
Expanding $D_1(\hat C)$ gives
\[
D_1(\hat C)=(Y_0-Q_p\hat C)^*(Y_0-Q_p\hat C),
\]
or, equivalently,
\[
D_1(\hat C)=Y_0^*Y_0-2Y_0^*Q_p\hat C+(\hat C)^*Q_p^*Q_p\hat C,
\]
because $Y_0^*Q_p\hat C=(\hat C)^*Q_p^*Y_0$. Furthermore
\[
\frac{\partial}{\partial \hat C }D_1(\hat C)=-2Q_p^*Y_0+2Q_p^*Q_p\hat C.
\]
Setting the derivative to zero, $\frac{\partial}{\partial \hat C }D_1(\hat C)=0$, we get
\[
Q_p^*Q_p\hat C=Q_p^*Y_0.
\]
Since $rank Q_p=p$, the matrix $Q_p^*Q_p$ is invertible and the solution is given by
\[
\hat C=(Q_p^*Q_p)^{-1}Q_p^*Y_0.
\]
The proof is completed.

\subsection*{Numerical example}
\textsl{Example 3.1.}\\\\
Assume the system \eqref{eq1} with non-consistent initial conditions of type \eqref{eq2} and $V_k=0_{5,1}$. Let
\[
F=\left[\begin{array}{ccccc} 0&0&-1&0&1\\0&-1&-1&1&1\\-1&-1&1&1&0\\0&1&2&0&-2\\0&0&0&0&0\end{array}\right],
\]
\[
G=\frac{1}{5}\left[\begin{array}{ccccc} -5&0&8&5&-3\\-11&-1&14&11&-8\\-2&-2&2&2&0\\11&2&-14&-11&8\\-5&0&10&5&-5\end{array}\right]
\]
and
\[
Y_0=\left[\begin{array}{c}1\\1\\1\\1\\1\end{array}\right].
\]
Then det($sF-G$)=$s(s-\frac{1}{5})(s-\frac{2}{5})$ and the pencil is regular. The three finite eigenvalues ($p$=3) of the pencil are 0, $\frac{1}{5}$, $\frac{2}{5}$ and the Jordan matrix $J_p$ has the form
\[
J_p=\left[\begin{array}{ccc}0&0&0\\0&\frac{1}{5}&0\\0&0&\frac{2}{5}\end{array}\right].
\]
By calculating the eigenvectors of the finite eigenvalues we get the matrix
\[
Q_p=\left[\begin{array}{ccc}1&0&0\\0&1&1\\0&0&0\\1&0&1\\0&0&1\end{array}\right].
\]
From Proposition 3.1, the initial conditions are non-consistent and thus from Theorem 3.1 and \label{eq7}, the optimal solution of the initial value problem \eqref{eq1} -- \eqref{eq2} is given by
\[
\hat Y_k=Q_pJ_p^k(Q_p^TQ_p)^{-1}Q_p^TY_0.
\]
Where 
\[
(Q_p^TQ_p)^{-1}=\frac{1}{3}\left[\begin{array}{ccc} 2&1&-1\\1&5&-2\\-1&2&-2\end{array}\right],
\]
\[
Q_p^TY_0=\left[\begin{array}{c} 2\\1\\3\end{array}\right]
\]
and
\[
(Q_p^TQ_p)^{-1}Q_p^TY_0=\frac{1}{3}\left[\begin{array}{ccc} 2\\1\\2\end{array}\right].
\]
Hence
\[
\hat Y_k=\frac{1}{3\cdot 5^k}\left[\begin{array}{ccc} 0\\3\\0\\4\\4\end{array}\right].
\]
\textsl{Example 3.2.}
\\\\
Assume the system \eqref{eq1} with initial conditions of type \eqref{eq2}. Let
\[
F=\left[\begin{array}{cc} 1&1\\1&1\end{array}\right],
\]
\[
G=\frac{1}{5}\left[\begin{array}{cc} 1&-2\\-2&0\end{array}\right]
\]
and
\[
Y_0=\left[\begin{array}{c}2\\3\end{array}\right].
\]
Then det($sF-G$)=$s+\frac{4}{5}$ and the pencil is regular. The finite eigenvalue ($p$=1) of the pencil is $-\frac{4}{5}$ and the Jordan matrix $J_p$ has the form
\[
J_p=-\frac{4}{5}.
\]
By calculating the eigenvectors of the finite eigenvalues we get the matrix
\[
Q_p=\left[\begin{array}{c}2\\3\end{array}\right].
\]
It is easy to observe that 
\[
Y_0\in colspan Q_p.
\]
From Proposition 3.1, the initial conditions are consistent and thus from \eqref{eq6} the unique solution of the initial value problem \eqref{eq1} -- \eqref{eq2} is given by 
\[
Y_k=Q_pJ_p^kC.
\]
or, equivalently,
\[
Y_k=-\frac{4^k}{5^k}\left[\begin{array}{c}2\\3\end{array}\right]C.
\]
Where $C$ is the unique solution of the linear system 
\[
\left[\begin{array}{c}2\\3\end{array}\right]=\left[\begin{array}{ccc}2\\3\end{array}\right]C
\]
and thus
\[
C=1,
\]
i.e.
\[
Y_k=-\frac{4^k}{5^k}\left[\begin{array}{c}2\\3\end{array}\right].
\]
\textsl{Example 3.3.}
\\\\
We assume the system \eqref{eq1} as in example 3.2 but with different initial conditions. Let
\[
Y_0=\left[\begin{array}{c}2.00001\\2.99999\end{array}\right].
\]
It is easy to observe that 
\[
Y_0\notin colspan Q_p.
\]
This means that from Proposition 3.1 the initial conditions are non-consistent and thus from Theorem 3.1 and \eqref{eq7} an optimal solution of the initial value problem is
\[
\hat Y_k=Q_pJ_p^k(Q_p^TQ_p)^{-1}Q_p^TY_0.
\]
Where 
\[
(Q_p^TQ_p)^{-1}=\frac{1}{13},
\]
\[
Q_p^TY_0=12.99999
\]
and
\[
(Q_p^TQ_p)^{-1}Q_p^TY_0=\frac{12.99999}{13}.
\]
Hence
\[
\hat Y_k=-\frac{12.9999\cdot4^k}{13\cdot5^k}\left[\begin{array}{c}2\\3\end{array}\right].
\]
It is worth to observe that the above optimal solution of the system \eqref{eq1} with non-consistent initial conditions $\left[\begin{array}{c}2.00001\\2.99999\end{array}\right]$, is almost equal to the unique solution of the same system with consistent initial conditions $\left[\begin{array}{c}2\\3\end{array}\right]$. This is because the perturbation to the non-consistent initial conditions \eqref{eq2} accordingly $min\left\|Y_0-\hat Y_0\right\|_2$=0.00018 is very small.

\section*{Conclusions}
In this article we focused and provided properties for optimal solutions of  a linear generalized discrete time system whose coefficients are square constant matrices and its pencil is regular. 


\end{document}